\documentclass[12pt]{article}
\usepackage{tcolorbox}
\usepackage{amssymb,amsmath,bm}
\usepackage[colorlinks=true, urlcolor=blue,linkcolor=blue, citecolor=blue]{hyperref}
\usepackage{fancyhdr}
\usepackage[top=2cm]{geometry}
\usepackage{mathrsfs,verbatim}
\usepackage{caption,cleveref }
\usepackage{graphicx, float}
\usepackage{constants,color,appendix,  mathtools, tikz, xcolor}
\usepackage{enumerate}
\usepackage{bbm}
\usepackage{tkz-euclide}
\usepackage{tikz}
\usepackage{longtable}

\usetikzlibrary{trees}
\usepackage{pdfpages, hhline}

\usetikzlibrary{calc}
\usepackage{lipsum}

\newcommand{\E}{\mathbb{E}}
\newcommand{\Pro}{\mathbb{P}}

\usepackage{appendix}
\begin{document}
\newcommand{\bea}{\begin{eqnarray}}
\newcommand{\ena}{\end{eqnarray}}
\newcommand{\beas}{\begin{eqnarray*}}
\newcommand{\enas}{\end{eqnarray*}}
\newcommand{\beq}{\begin{equation}}
\newcommand{\enq}{\end{equation}}
\def\qed{\hfill \mbox{\rule{0.5em}{0.5em}}}
\newcommand{\bbox}{\hfill $\Box$}
\newcommand{\ignore}[1]{}
\newcommand{\ignorex}[1]{#1}
\newcommand{\wtilde}[1]{\widetilde{#1}}
\newcommand{\mq}[1]{\mbox{#1}\quad}
\newcommand{\bs}[1]{\boldsymbol{#1}}
\newcommand{\qmq}[1]{\quad\mbox{#1}\quad}
\newcommand{\qm}[1]{\quad\mbox{#1}}
\newcommand{\nn}{\nonumber}
\newcommand{\Bvert}{\left\vert\vphantom{\frac{1}{1}}\right.}
\newcommand{\To}{\rightarrow}
\newcommand{\supp}{\mbox{supp}}
\newcommand{\law}{{\cal L}}
\newtheorem{problem}{Problem}[section]
\newtheorem{exercise}{Exercise}[section]
\newtheorem{theorem}{Theorem}[section]
\newtheorem{eg}{Example}[section]
\newtheorem{thm}{Theorem}[section]
\newtheorem{lem}{Lemma}[section]
\newtheorem{soln}{Solution}[section]
\newtheorem{propn}{Proposition}[section]
\newtheorem{ex}{Exercise}[section]
\newtheorem{corollary}{Corollary}[section]
\newtheorem{conjecture}{Conjecture}[section]
\newtheorem{proposition}{Proposition}[section]
\newtheorem{lemma}{Lemma}[section]
\newtheorem{definition}{Definition}[section]
\newtheorem{example}{Example}[section]
\newtheorem{remark}{Remark}[section]
\newtheorem{solution}{Solution}[section]
\newtheorem{case}{Case}[section]
\newtheorem{condition}{Condition}[section]
\newcommand{\pf}{\noindent {\bf Proof:} }
\newcommand{\proof}{\noindent {\textbf{Proof:}} }
\frenchspacing

\title{\bf  On a cover time  problem on a dynamic graph with steps at random times}

\author{Yunus Emre Demirci\footnote{Queen's University, Department of Mathematics and Statistics, Kingston, ON, Canada, email:  21yed@queensu.ca} \hspace{0.2in} \"{U}m\.{i}t I\c{s}lak\footnote{Bo\u{g}azi\c{c}i University, Department of Mathematics, Istanbul, Turkey, email: umit.islak1@boun.edu.tr}  \hspace{0.2in} Mehmet Akif Y{\i}ld{\i}z\footnote{‡University of Amsterdam, Korteweg-de Vries Institute for Mathematics, Amsterdam, the Netherlands,
e-mail: m.a.yildiz@uva.nl}} 
\vspace{0.25in}

\maketitle

\begin{abstract}
We introduce a cover time problem for   random walks on   dynamic graphs in which the graph expands in time and the walker moves at   random times. Number of nodes covered at a certain time and number of  returns to original states are analyzed in the resulting model.

\bigskip

Keywords:  Dynamic graph, random walk, cover time, coupon collecting. 


\end{abstract}

\tikzstyle{level 1}=[level distance=2.75cm, sibling distance=5.65cm]
\tikzstyle{level 2}=[level distance=3cm, sibling distance=2.75cm]
\tikzstyle{level 3}=[level distance=3.9cm, sibling distance=1.5cm]

\tikzstyle{bag} = [text width=10em, text centered] 
\tikzstyle{end} = [circle, minimum width=3pt,fill, inner sep=0pt]


\section{Introduction}

Let us begin by recalling the well-known  coupon collector problem.
Consider  a set containing $m$ distinct objects (``coupons"),
for example, pictures of soccer players. The collector samples from the set with
replacement. On each trial she has a fixed probability $p_i$  of drawing a type  $i$ object, independently of all past events. A well-known question is  then the expected  time to collect all coupons, though various other probabilities/expectations  can be of interest. The coupon collecting emerges in distinct sciences, and    has various applications. See \cite{fs} for an elementary treatment of this problem. \cite{fgt} is a classical paper on the subject with relations  to other well known problems in probability theory and computer science.  

In our setting, the coupon collector problem  can be considered as a special case of the cover time problem for random walks on graphs (which is the complete graph $K_n$ in  case of  the  coupon collector problem).  
 In this general setup, one is 
interested in the cover time required for all vertices of a graph to be  visited at least once when the walker   traverses randomly on the vertices according 
to certain probabilities.    The cover time in random walks (which is a  Markov chain) has been investigated in
numerous studies in the literature.
 Two  general references for Markov chains are   \cite{MT} and  \cite{norris}. The reader can find various  results on cover times   in  \cite{aldousfill} and \cite{LPW}. The former also contains various related conjectures. 
The standard cover time problem has several  variations in literature such as the edge cover problem \cite{zuckerman},  having multiple walkers \cite{broder} or covering trees \cite{brightwell}.

Our purpose here is to give yet another  very simple variation of the classical cover time problem on graphs. On contrary to the classical setting, we are interested in a dynamic version in which the walker moves in a graph which expands in time.   Such dynamic graphs arise in various areas, for example, in collaboration networks,  acquaintance graphs and the web graph. Our main motivation was indeed the   web graph in which the vertices correspond to web sites, and the edges correspond to links between them. This is clearly  a dynamic graph since new sites open each day and they link to previous existing pages. Our setup below is that a random walker starts at a certain web page initially and wanders around the evolving graph, and we try to understand statistics such as the number of nodes visited until a specific time. We are curious about the possible uses of such modeling for a search engines which has the goal to discover as many web pages as possible by doing a random walk. The  approaches in references   \cite{evolvinggraph} and  \cite{evolvingweb}  are    same direction, but both  are different than our model. Also, after posting our paper on arxiv, we were got in touch by Dr. Videla about a recent work following the same model we discuss below \cite{videla}. Due to this,  there had been some minor changes in our manuscript, with references to the cited paper. 
The model we discuss below is as follows:  
\begin{itemize}
	\item[(i)]  Assume there are $k_0$ (call $v_1$, $v_2, \ldots, v_{k_0}$) vertices at the beginning (at time $0$) for some $k_0\geq 3$, and that the initial configuration is a complete graph on these $k_0$ vertices. We also assume that there is a walker at one of these vertices.
	\item[(ii)] For each $t\in\mathbb{N}$, create a new vertex $v_{k_0+t}$ at time $t$, and join it to all present vertices.
	\item[(iii)] The walker moves   to one of the remaining current vertices at random where the time intervals between moves are exponentially distributed with parameter $1 / \lambda$.
\end{itemize}
Taking this construction as our base ground, we study mainly two statistics throughout the paper: (1) Number of nodes visited up to a given time, (2) Number of visits to  a specified initial node up to a certain time.  
Note that the arrival of a new vertex  could  also be modeled as an independent Poisson process, but we decided to keep  the deterministic one  in our   analysis because otherwise  the notations turn out to be  cumbersome.  In either case as a technical remark, let us note that since the graph structure evolves in time, the underlying  probability  space changes accordingly. For relevant  discussions of the construction of the process, the reader may check \cite{videla}.

As already noted, our model is a special case of the one proposed by \cite{videla}, but below we give a more detailed  analysis of the   statistics  we study due to this specialization. Before the main discussion, let us also briefly tell about the references \cite{evolvinggraph} and  \cite{evolvingweb} where   cover time problems for   dynamic graphs are analyzed. In   \cite{evolvingweb}, the authors consider  two similar models for growing the graph with an  inspiration from web graphs. In both models, there is a sequence of connected graphs $\mathcal{G}=\{G(t):t=1,2,\ldots\}$ for which $G(t)$ is constructed from $G(t-1)$ by including a new vertex and a fixed number of edges between the new vertex and the vertices of $G(t-1)$. Here, the initial graph $G(1)$ consists of a single vertex with some self-loops.   While the neighbors of the new vertex is chosen independently and uniformly with the possibility of multiple connections in the first model, the new vertex is linked with previous vertices with probability proportional to their degrees in the second one. For both models, there is a walker traversing on the vertices where she walks for a fixed length at each step. 
Following these models, the authors examine the expectation of the number of vertices which have not been visited by the walker up to a certain step. In particular, the authors study certain asymptotics related to the portion of the vertices   not   visited by the walker.

On the other hand,  another dynamic model   was examined in \cite{evolvinggraph} which is much different than our model  since in their setting the graphs evolves in time but on a fixed set of vertices.   It is known that the covering of static graphs always requires a polynomial time. In their work in \cite{evolvinggraph}, the authors prove that this is not the case for evolving graph models. They construct a sequence of stars in a way that the hitting time between two specific vertices is exponential in terms of the size of stars, which implies that the cover time is not polynomial.


  
   The rest of the paper is organized as follows. Next section presents our foundings on the number of   nodes covered  at a certain time $T$ which we call $N_T$. We study the expectation, the variance and certain asymptotics related to this statistic. In Section \ref{sec:numvis} we will be looking at probabilities and expectations related to the number of visits to the initial vertices. 

\section{Number of   nodes covered  at a certain time}\label{sec:N-T}

Unless otherwise stated, we consider the model introduced in Introduction. 
For  $T \in \mathbb{N}$, let $N_T $ be the number of vertices visited at least once by the walker up to time $T$.  The purpose of this section is to analyze $\mathbb{E}[N_T]$ and $\mathrm{Var}(N_T)$ along with the asymptotics related to $N_T$.  
Below $\mathcal{H}(n)$ denotes finite sums of  harmonic series, i.e. $\mathcal{H}(n) = \sum_{i=1}^n \frac{1}{i}$.

\subsection{Expectation}

We begin with the analysis of $\mathbb{E}[N_T]$. Let us note that Theorem \ref{thm:exp} is  available in \cite{videla}. We include the details of our elementary treatment for the sake of complements.  Also, the argument given here will be useful later in Section \ref{sec:varc}.

\begin{theorem}\label{thm:exp} (i) We have 
\begin{equation}\label{eqn:exp}
	\E[N_T]=(k_0+T-1)-e^{-\lambda\cdot \mathcal{H}(k_0+T-2)}\Big((k_0-1)e^{\lambda\cdot \mathcal{H}(k_0-2)}+\sum_{j=k_0+1}^{k_0+T-1}e^{\lambda\cdot \mathcal{H}(j-2)}\Big).
\end{equation}
(ii) As $T \rightarrow \infty$, $\frac{\mathbb{E}[N_T]}{T} \sim \frac{\lambda}{1 + \lambda}.$
\end{theorem}

\textbf{Proof.} (i)
Let $T > 0$ and $u_j = u_j(T)$ be the probability that the vertex $j$ is visited at least once from the walker up to time $T$. Without loss of generality, assume we start at the vertex $v_1$. Clearly, $u_1=1$ and $u_s=0$ for $s\geq k_0+T$. Observe that 
 $
	\E[N_T]=  1+\sum_{j=2}^{k_0+T-1} u_j.
$ 

Let us now calculate the probability $u_j$ for   $1\leq j\leq T$. Let $M_t$ be the number of moves by the walker  on the time interval $[t-1,t)$ for $t\geq1$. We know $M_t$ has Poisson distribution with parameter $\lambda$. Also, it is well-known that $M_t$ and $M_s$ are independent for $t\neq s$. Therefore, for any $1\leq j\leq T$ and $m_j,m_{j+1},\ldots,m_T\in\mathbb{N}_{\geq0}$, we get  $\Pro(M_j=m_j,M_{j+1}=m_{j+1},\ldots,M_T=m_t)=\prod_{s=j}^{T}\dfrac{e^{-\lambda}\lambda^{m_s}}{m_s!}.$ 

Now, suppose $M_1=m_1$, $M_2=m_2,\ldots,M_T=m_T$ are given where $m_1,m_2,\ldots,m_T\in\mathbb{N}_{\geq0}$. If $2\leq j\leq k_0$, then the walker does not visit the vertex $v_j$ up to time $T$ with probability $\Big(\dfrac{k_0-2}{k_0-1}\Big)^{m_1}\cdot\Big(\dfrac{k_0-1}{k_0}\Big)^{m_2}\cdot\Big(\dfrac{k_0}{k_0+1}\Big)^{m_3}\cdots \Big(\dfrac{k_0+T-3}{k_0+T-2}\Big)^{m_T}.$
Let us call this quantity as $R[j;m_1,m_2,\ldots,m_T]$. Similarly, if $k_0+1\leq j\leq k_0+T-1$, then the walker does not visit the vertex $v_j$ up to time $T$ with probability $\Big(\dfrac{j-2}{j-1}\Big)^{m_{j-k_0+1}}\cdot\Big(\dfrac{j-1}{j}\Big)^{m_{j-k_0+2}}\cdot\Big(\dfrac{j}{j+1}\Big)^{m_{j-k_0+3}}\cdots \Big(\dfrac{k_0+T-3}{k_0+T-2}\Big)^{m_T}.$
We  write $S[j;m_1,m_2,\ldots,m_T]$ for this last probability.

Now, for $2\leq j\leq k_0$, by   conditioning on $M_1$, $M_2,\ldots,M_T$, we have
$$	1-u_j = \sum_{m_1,m_2,\ldots,m_T\in\mathbb{N}_{\geq0}} R[j;m_1,m_2,\ldots,m_T]\cdot \Pro\big(M_i=m_i \text{ for }1\leq i \leq T\big)
$$
On the other hand, observe that we can express $R[j;m_1,m_2,\ldots,m_T]$ as $ \prod_{s=1}^{T} \big(1 - \frac{1}{k_0+s-2}\big)^{m_s}$. Therefore, we can write 
\begin{eqnarray*}
	1-u_j &=& \sum_{m_1,m_2,\ldots,m_T\in\mathbb{N}_{\geq0}} \Bigg(\prod_{s=1}^{T}\big(1 - \frac{1}{k_0+s-2}\big)^{m_s}\Bigg) \cdot \dfrac{e^{-\lambda}\lambda^{m_s}}{m_s!} \\
	&=&\sum_{m_1,m_2,\ldots,m_T\in\mathbb{N}_{\geq0}} \prod_{s=1}^{T} e^{-\lambda} \dfrac{\big(\lambda - \frac{\lambda}{k_0+s-2}\big)^{m_s}}{m_s!} \\
	&=& \sum_{m_1,m_2,\ldots,m_T\in\mathbb{N}_{\geq0}} e^{-\lambda T}\cdot\prod_{s=1}^{T}  \dfrac{\big(\lambda - \frac{\lambda}{k_0+s-2}\big)^{m_s}}{m_s!} 
	=  e^{-\lambda T}\cdot\sum_{m_1,m_2,\ldots,m_T\in\mathbb{N}_{\geq0}}\prod_{s=1}^{T}  \dfrac{\big(\lambda - \frac{\lambda}{k_0+s-2}\big)^{m_s}}{m_s!} \\
	&=&  e^{-\lambda T}\cdot\prod_{s=1}^{T} \sum_{l=0}^{\infty}  \dfrac{\big(\lambda - \frac{\lambda}{k_0+s-2}\big)^{l}}{l!} 
	 =  e^{-\lambda T} \cdot  \prod_{s=1}^{T} e^{\lambda - \frac{\lambda}{k_0+s-2}}= e^{-\lambda\cdot\big(\frac{1}{k_0-1}+\frac{1}{k_0}+ \cdots +\frac{1}{k_0+T-2}\big)}
\end{eqnarray*} 

In a similar way, for $k_0+1\leq j\leq k_0+T-1$, again   conditioning  on $M_{k-j_0+1}, M_{k-j_0+2}, \ldots, M_T$ gives 
$$
1-u_j = \sum_{m_{j-k_0+1},m_{j-k_0+2},\ldots,m_T\in\mathbb{N}_{\geq0}} S[j;m_1,m_2,\ldots,m_T]\cdot \Pro\big(M_i=m_i \text{ for }j-k_0+1\leq i \leq T\big).$$
Noting that we can express $S[j;m_1,m_2,\ldots,m_T]$ as $\prod_{s=j-k_0+1}^{T}\big(1 - \frac{1}{k_0+s-2}\big)^{m_s}$,  computations as in previous case then gives 
$	1-u_j = e^{-\lambda\cdot\big(\frac{1}{j-1}+\frac{1}{j}+\cdots +\frac{1}{k_0+T-2}\big)}, \, k_0+1\leq j\leq k_0+T-1.$

As a result,  
\begin{equation}\label{def:uj}
u_j=\begin{cases}
1-e^{-\lambda\cdot\big(\mathcal{H}(k_0+T-2)-\mathcal{H}(k_0-2)\big)}, & \text{ if } 2\leq j\leq k_0. \\
1-e^{-\lambda\cdot\big(\mathcal{H}(k_0+T-2)-\mathcal{H}(j-2)\big)}, & \text{ if }  k_0+1\leq j\leq k_0+T-1.   
\end{cases}
\end{equation}

\noindent Therefore,  recalling $\E[N_T]=1+\displaystyle{\sum_{j=2}^{k_0+T-1} u_j}$, the result follows.

(ii) The asymptotic assertion follows from the inequalities $ \frac{1}{2(m+1)}  + \ln m + \gamma < \mathcal{H}(m) <  \frac{1}{2m}  + \ln m + \gamma$, $m \in \mathbb{N}, $ where $\gamma = 0.5772\cdots$ is the Euler-Mascheroni constant. 
	 \hfill $\square$

\subsection{Variance}\label{sec:varc}
 
\begin{theorem}\label{thm:events-A-j} (i) For  $2 \leq r<s\leq k_0+T-1$, the events $A_r$ and $A_s$ 
are independent. 

(ii) We have 
$\mathrm{Var}(N_T)=\displaystyle{\sum_{j=2}^{k_0+T-1} (1-u_j)\cdot u_j },$
 where $u_j$'s are as given in \eqref{def:uj}.

(iii) We have $\mathrm{Var}(N_T) \sim \frac{\lambda}{(\lambda+1) (2 \lambda + 1)} T$.
\end{theorem}

\textbf{Proof.} (i)-(ii)
  Let $A_j = A_j(T)$ be the event that the vertex $j$ is not visited by the walker up to time $T$. Suppose $M_1=m_1, M_2=m_2,\ldots,M_T=m_T$ are given where $m_1,m_2,\ldots,m_T\in\mathbb{N}_{\geq0}$.
Let us assume that $k_0+1\leq r<s\leq k_0+T-1$. Then, the event $A_r\cap A_s$ occurs with probability
$\Big(\dfrac{r-2}{r-1}\Big)^{m_{r+1}}\Big(\dfrac{r-1}{r}\Big)^{m_{r+2}}\cdots\Big(\dfrac{s-3}{s-2}\Big)^{m_{s}}\Big(\dfrac{s-3}{s-1}\Big)^{m_{s+1}}\Big(\dfrac{s-2}{s}\Big)^{m_{s+2}}\cdots\Big(\dfrac{k_0+T-4}{k_0+T-2}\Big)^{m_{T}}.$ \\
Thus, with a    conditioning argument similar to the expectation case, we may write
$$
\Pro(A_r\cap A_s)=e^{-\lambda\cdot\big(\frac{1}{r-1}+\frac{1}{r}+\cdots\frac{1}{s-2}+\frac{2}{s-1}+\frac{2}{s}+\cdots+\frac{2}{k_0+T-2}\big)}
=e^{-\lambda\cdot\big(\frac{1}{r-1}+\frac{1}{r}+\cdots\frac{1}{k_0+T-2}\big)}  e^{-\lambda\cdot\big(\frac{1}{s-1}+\frac{1}{s}+\cdots + \frac{1}{k_0+T-2}\big)}.
$$
As a result, we have $\Pro(A_r\cap A_s)=(1-u_r)(1-u_s)=\Pro(A_r)\Pro(A_s),$ for $k_0+1\leq r<s\leq k_0+T-1$. In a similar   way, it can be seen that the equality $\Pro(A_r\cap A_s)=\Pro(A_r)\Pro(A_s)$ holds for all $2\leq r<s\leq k_0+T-1$.

Define $U_T$ as the number of unvisited vertices among $\{v_1,v_2,\ldots,v_{k_0+T-1}\}$ by the walker up to time $T$. Since $N_T+U_T=k_0+T-1$, we have $\mathrm{Var}(N_T)=\mathrm{Var}(U_T)$. Now, by using the equalities $U_T=\displaystyle{\sum_{j=2}^{k_0+T-1}\mathbf{1}(A_j)}$, $\Pro(A_j)=1-u_j$, and the fact that equality $\Pro(A_r\cap A_s)=\Pro(A_r)\Pro(A_s)$ holds for all $2\leq r<s\leq k_0+T-1$, we get 
$\mathrm{Var}(U_T) = \displaystyle{\sum_{j=2}^{k_0+T-1} (1-u_j)u_j}$, from which   $\mathrm{Var}(N_T)=\displaystyle{\sum_{j=2}^{k_0+T-1} (1-u_j)\cdot u_j }$ follows.  

(iii) Again follows from $ \frac{1}{2(m+1)}  + \ln m + \gamma < \mathcal{H}(m) <  \frac{1}{2m}  + \ln m + \gamma, \, m \in \mathbb{N}$.
\hfill $\square$


 	
\begin{remark} (Cover time of a fixed subset) The independence in previous result can also be used for studying the classical cover time of a fixed subset of the vertices. Let us briefly include a discussion about this. Letting $k\geq k_0$ be fixed, we study  the cover time of the set $\{1,2,\ldots,k\}$.  Let $C_k$ be the first time at which all the vertices in $\{1,2,\ldots,k\}$ have been visited. For  convenience, we assume all the moves by the walker on the time interval $(t-1,t]$ are done at time $t$. By  assumption, $C_k$ can take only integer values. Let $T\geq k-k_0+1$, and recall that the events $B_j(T) : = \{\text{vertex } j \text{ is visited by time } T \}$ are independent by Theorem~\ref{thm:events-A-j} for $2\leq r<s\leq k_0+T-1$. Hence, 
	\begin{align*}
	\Pro(C_k\geq T)
	=1-\Pro\left(B_2(T)\cap B_3(T)\cap \cdots \cap B_{k}(T)\right)
	=1-u_2(T)u_3(T)\cdots u_k(T),
	\end{align*} 
where $u_i(T)$ are given as in \eqref{def:uj} (note the dependence on $T$). Then, the cover time is
\begin{align*}
\mathbb{E}[C_k] = \sum_{j =1}^{\infty} \mathbb{P}(C_k \geq j) &= \sum_{j=1}^{k - k_0}  \mathbb{P}(C_k \geq j)  + \sum_{T = k - k_0+1}^{\infty}  \mathbb{P}(C_k \geq T)  \\
&= (k  - k_0)+ \sum_{T=k-k_0+1}^{\infty}\left(1-\prod_{j=2}^{k}u_j(T)\right).
\end{align*}

Now, specializing to the case $k=k_0$, and noting   $u_2(T)=\cdots=u_{k_0}(T)$ for any $T$, this can be rewritten as, 
\begin{align*}
\E[C_k]&=\sum_{T=1}^{\infty}\left(1-u_2(T)^{k_0-1}\right)=\sum_{T=1}^{\infty}\left(1-\left( 1-e^{-\lambda\cdot\left(\mathcal{H}(k_0+T-2)-\mathcal{H}(k_0-2)\right)} \right)^{k_0-1}\right)\\
&=\sum_{T=1}^{\infty}\sum_{j=1}^{k_0-1}\binom{k_0-1}{j}(-1)^{j+1}e^{-\lambda j\cdot\left(\mathcal{H}(k_0+T-2)-\mathcal{H}(k_0-2)\right)}\\
&=\sum_{j=1}^{k_0-1}\binom{k_0-1}{j}(-1)^{j+1}\sum_{T=1}^{\infty}e^{-\lambda j\cdot\left(\mathcal{H}(k_0+T-2)-\mathcal{H}(k_0-2)\right)}.
\end{align*}
\end{remark}





 
\subsection{Asymptotics for $N_T$}

The purpose of this section is to provide  a  law of large numbers  and  a central limit  theorem for $N_T$. Towards the former,  we first derive a deviation inequality for $N_T$. For this purpose recall that the events $B_j : = \{\text{vertex } j \text{ is visited by time } T \}$ are independent. Letting then $Y_j$ be the indicator of $B_j$, and using Azuma-Hoeffding   inequality, for any $t > 0$,  we have 
$$
\mathbb{P}(|N_T - \mathbb{E}[N_T]| \geq t) = \mathbb{P} \left(\left|\sum_{j=1}^{k_0+T} Y_j - \mathbb{E}[N_T] \right| \geq t \right) 
\leq \exp \left(- \frac{2t^2}{\sum_{j=1}^{k_0+T} 1 } \right) 
= \exp \left(- \frac{2t^2}{k_0+T} \right).
$$
In other words, letting $\alpha_T = \mathbb{E}[N_T / T]$, for any $\epsilon > 0$, we have $$\mathbb{P}\left(\left|\frac{N_T}{T} - \alpha_T \right| \geq \epsilon \right) \leq \exp \left( - \frac{2 \epsilon^2 T^2}{k_0 + T} \right), \quad T > 0.$$   Then, for any $\epsilon > 0$ $$\sum_{T= 1}^{\infty}\mathbb{P}\left(\left|\frac{N_T}{T} - \alpha_T \right| \geq \epsilon \right) \leq \sum_{T= 1}^{\infty} \exp \left( - \frac{2 \epsilon^2 T^2}{k_0 + T} \right) < \infty.$$ Therefore, the first Borel-Cantelli lemma and Theorem \ref{thm:exp}(i) yields 
\begin{theorem}
$\frac{N_T}{T} \rightarrow_{a.s.}  \frac{\lambda}{1 + \lambda} $. 
\end{theorem}
 
Next we  prove a central limit theorem for $N_T$. For this purpose we will be using the following standard result.  

\begin{theorem} 
Suppose $X_1,X_2,\ldots$ be a sequence of mean zero, independent  random variables such that $\mathbb{E}[X_i^4] < \infty$ for each $i$.  If $W_n = \frac{ \sum_{i=1}^n X_i}{\sigma}$ and  $\sigma^2 = \mathrm{Var} \left(\sum_{i=1}^n X_i \right)$, then $$d_W (W_n, \mathcal{G}) \leq \frac{1}{\sigma^3 } \sum_{i=1}^n \mathbb{E}|X_i|^3 + \frac{C \sqrt{ \sum_{i=1}^n \mathbb{E}[X_i^4]}}{\sigma^2}, $$ where $C$ is some constant independent of $n$, $d_W$ is the Wasserstein distance between probability measures and $\mathcal{G}$ is a standard Gaussian random variable.
\end{theorem}

In our case,  we set $X_i = \mathbf{1} (B_i ) - \mathbb{P} (B_i)$, $i \geq 1$,  where  $B_j : = \{\text{vertex } j \text{ is visited by time } T \}$. Then the following central limit theorem follows after straightforward manipulations  with keeping 
$\mathrm{Var}(N_T) \sim \frac{\lambda}{(\lambda+1) (2 \lambda + 1)} T$ in mind. 

\begin{theorem}
We have 
$d_W \left(\frac{N_T - \mathbb{E}[N_T]}{\mathrm{Var} (N_T)}, \mathcal{G} \right) \leq \frac{C}{\sqrt{n}},$ for some constant $C$ independent of $n$. In particular, $\frac{N_T - \mathbb{E}[N_T]}{\mathrm{Var} (N_T)}$ converges in distribution to a standard normal random variable. 
\end{theorem}
 
 \begin{remark}
 To keep the discussion simple, we take the node generation times   deterministic in this manuscript. However, some elementary manipulations could be used to carry the results to random generation times. A related brief discussion can be found in the arxiv version of the paper. 
 \end{remark}

 \section{Number of visits to the initial nodes}\label{sec:numvis}
 
The main purpose of this section  to study the expected number of visits to the initial state $v_1$ up to time $T$ when there are initially $k_0 > 0$ many vertices. This is Theorem \ref{thm:mainexp} below.   Throughout the way we  will compute various other related probabilities. 
First, we begin with the  calculation of  the probability of not visiting the initial vertex a second time until a specified time.

When there are $k_0$ vertices at the beginning, let us define $P(T,k_0)$ to be the probability of not visiting the initial vertex a second time up to time $T$. Obviously $P(0, k_0) = 1$.

\begin{proposition}
For all $T > 1$,  $P(T,k_0) = e^{-\lambda  T}+\sum_{i=0}^{T-1} \left(1+ \frac{1}{k_0-2+i} \right) e^{-\lambda \left(i+\sum_{j=1}^{T-i}\frac{1}{k_0-2+i+j}\right)}$.
  \end{proposition}

\textbf{Proof.}
Recall that $M_t$ is  the number of moves on the interval $[t-1,t)$ for $t\geq1$, and that after a unit of time, a new node is created. If $M_1=0$, a unit time later the walker  is still on the initial vertex and the number of vertices is $k_0+1$. If $M_1>0$, a unit time later it should not have visited initial vertex a second time and the number of vertices is $k_0+1$. 

We   have 
\begin{eqnarray}\label{ptk0eq}
P(T,k_0) &=& P(M_1=0) P(T-1,k_0+1) \\
\nonumber &+&P(M_1\geq 1 \text{ and no visits to $v_1$ in $M_1$ steps}) (1-u_*),
\end{eqnarray}
where $u_*$ is the probability that we do have a revisit to the state where we are at time $1$ during the interval $(1,T]$. 
Here is a bit more explanation for the reasoning here. The first term on the right-hand side of \eqref{ptk0eq} is due to  the case where the walker did not have any moves during the initial unit interval.  Regarding the second term on the right-hand side, the walker had a move  in the first unit time interval and ended up  her move at some other vertex. This ``other" vertex is now considered to be the new initial vertex and we then look at the  probability of not visiting it for which we write   $1-u_*$. 

Continuing our computation, we then have 
\begin{eqnarray*}
P(T,k_0) &=& e^{-\lambda} P(T-1,k_0+1) + \sum_{r=1}^{\infty} \frac{e^{-\lambda}\lambda^r}{r!} \left(\frac{k_0-2}{k_0-1}\right)^{r-1}  \left(e^{-\lambda \left(\frac{1}{k_0}+\cdots+\frac{1}{k_0+T-2}\right)}\right)\\
&=& e^{-\lambda} P(T-1,k_0+1) + 
e^{-\lambda} \sum_{r=1}^{\infty} \frac{\left(\lambda \frac{k_0-2}{k_0-1} \right)^r}{r!} \left(\frac{k_0-1}{k_0-2}\right)  \left(e^{-\lambda \left(\frac{1}{k_0}+\cdots+\frac{1}{k_0+T-2}\right)}\right)\\
&=& e^{-\lambda} P(T-1,k_0+1) + 
e^{-\lambda} e^{\lambda \frac{k_0-2}{k_0-1}} \left(\frac{k_0-1}{k_0-2}\right)  \left(e^{-\lambda \left(\frac{1}{k_0}+\cdots+\frac{1}{k_0+T-2}\right)}\right)\\
&=& e^{-\lambda} P(T-1,k_0+1) +\left(\frac{k_0-1}{k_0-2}\right)  \left(e^{-\lambda \left(\frac{1}{k_0-1}+\cdots+\frac{1}{k_0+T-2}\right)}\right)\\
\end{eqnarray*}
Iterating this recursion we obtain 
\begin{eqnarray*}
P(T,k_0)&=&\left(\frac{k_0-1}{k_0-2}\right)\left(e^{-\lambda \left(\frac{1}{k_0-1}+\cdots+\frac{1}{k_0+T-2}\right)}\right) \\
&& +  e^{-\lambda}\left(\left(\frac{k_0}{k_0-1}\right)\left(e^{-\lambda\left(\frac{1}{k_0}+\cdots+\frac{1}{k_0+T-2}\right)}\right) + e^{-\lambda} P(T-2,k_0+2)\right) \\ &=& \cdots  \\ 
&=& 
 \frac{k_0-1}{k_0-2}e^{-\lambda\left(\frac{1}{k_0-1}+\cdots+\frac{1}{k_0+T-2}\right)}+\frac{k_0}{k_0-1} e^{-\lambda\left(1+\frac{1}{k_0}+\cdots+\frac{1}{k_0+T-2}\right)}\\
 && + \frac{k_0+1}{k_0} e^{-\lambda\left(2+\frac{1}{k_0+1}+\cdots+\frac{1}{k_0+T-2}\right)}  + \cdots  +\frac{k_0+T-2}{k_0+T-3}e^{-\lambda\left(T-1+\frac{1}{k_0+T-2}\right)}+e^{-\lambda T}\\
&=&e^{-\lambda T}+\sum_{i=0}^{T-1} \frac{k_0-1+i}{k_0-2+i} e^{-\lambda\left(i+\sum_{j=1}^{T-i}\frac{1}{k_0-2+i+j}\right)}.
\end{eqnarray*}
 \hfill $\square$

Next  we look at the   probability of being at one of the initial vertices $v$ at the end of a certain amount of time  when the walker starts the walk from $v$ with some given probability.  

\begin{proposition}\label{propn:pr}
Assume that $p$ is the probability of starting the walk from vertex $v$ at the beginning. Let $p_r$ be the probability of being at vertex $v$ after $r$ steps of the random walk in a  unit time where no new vertex is generated. Then for all $r \geq 0$, we have
$p_r=\frac{1}{k_0}-\frac{1-pk_0}{k_0}\left(-\frac{1}{k_0-1}\right)^r$. 
\end{proposition}
\textbf{Proof.}
The result is clearly true for  $r=0$.
For the general case, in order   to be at the node $v$ after $r$ steps, it must be in some other node one step earlier and it must move to the node $v$ in the last step:
$p_r=(1-p_{r-1})\frac{1}{k_0-1}.$
Subtracting $p_{r-1}$ from $p_r$, this gives us for all $r\geq 2$
$$
p_r-p_{r-1}=(1-p_{r-1})\frac{1}{k_0-1}-(1-p_{r-2})\frac{1}{k_0-1}=
(p_{r-1}-p_{r-2})\left(-\frac{1}{k_0-1}\right).
$$
When we use this formula recursively, we   obtain that, 
$p_r-p_{r-1}=(p_{1}-p_{0})\left(-\frac{1}{k_0-1}\right)^{r-1}$ for all $r \geq 1$.
Summing up from 1 to $r$ gives
$$\sum_{i=1}^r(p_i-p_{i-1})=\sum_{i=1}^r \left((p_{1}-p_{0})\left(-\frac{1}{k_0-1}\right)^{i-1}\right), \; \, \text{ or, } \; \,
p_r=p_0+(p_1-p_0)\sum_{i=1}^r \left(-\frac{1}{k_0-1}\right)^{i-1}.
$$
We know that $p_1=(1-p_0)\frac{1}{k_0-1}$ and $p_0=p$, so 
$p_r=\frac{1}{k_0}-\frac{1-pk_0}{k_0}\left(-\frac{1}{k_0-1}\right)^r.$
\hfill $\square$

Now define $P_{k_0}(p)$ to be  the probability of being at vertex $v$ at the end of  the first unit time if the probability of starting from vertex $v$ equals to $p$ and the initial number of vertices equals to $k_0$.

\begin{proposition}
We have 
 $P_{k_0}(p)=\frac{1}{k_0}-\frac{1-pk_0}{k_0} e^{\frac{-\lambda k_0}{k_0-1}}.$

\end{proposition}

\textbf{Proof.}  Proof follows via the use of Proposition \ref{propn:pr} and the  following observations:
\begin{eqnarray*}
P_{k_0}(p) &=& \sum_{r=0}^\infty P(M_1=r) p_r 
 =  \sum_{r=0}^\infty \frac{e^{-\lambda}\lambda^r}{r!} p_r = \sum_{r=0}^\infty \frac{e^{-\lambda}\lambda^r}{r!} \left(\frac{1}{k_0}-\frac{1-pk_0}{k_0}\left(-\frac{1}{k_0-1}\right)^r\right) \\  &=&  \sum_{r=0}^\infty \frac{e^{-\lambda}\lambda^r}{r!} \frac{1}{k_0}-
\sum_{r=0}^\infty \frac{e^{-\lambda}\left(-\frac{\lambda}{k-1}\right)^r}{r!} \frac{1-pk_0}{k_0} = \frac{1}{k_0}-\frac{1-pk_0}{k_0} e^{\frac{-\lambda k_0}{k_0-1}}
\end{eqnarray*}
\hfill $\square$


Let next $Q(T,k_0)$ be the probability of being at the initial vertex at time $T$ when the number of vertices is $k_0$ at the beginning. 
There are two cases for it to be at vertex $v$ at time $T$. Either it will be at vertex $v$ at time $T-1$, and it will be still at vertex $v$ at time $T$. Or it will be in another vertex at time $T-1$ and jump to vertex $v$ during the last interval.

\begin{proposition}\label{propn:qtk}
We have
\begin{eqnarray*}
Q(T,k_0)=\left(\prod_{i=1}^{T}e^{-\lambda \frac{k_0+i-1}{k_0+i-2}}\right)
\left(1+\sum_{j=1}^{T}\left(\frac{\left(1-e^{-\lambda\frac{k_0+j-1}{k_0+j-2}}\right)}{(k_0+j-1)\prod_{i=1}^{j}e^{-\lambda \frac{k_0+i-1}{k_0+i-2}}}\right)\right).
\end{eqnarray*}

\end{proposition}

\textbf{Proof.}
Observe that $Q(0,k_0)=1$ and $Q(1,k_0)=P_{k_0}(1)=\frac{1}{k_0}+\frac{k_0-1}{k_0} e^{\frac{-\lambda k_0}{k_0-1}}.$
 If we recall the   observation just before the statement of the proposition, we see that for $T\geq1$,
\begin{small}
\begin{eqnarray*}
Q(T,k_0) & = & \mathbb{P}(\text{being at initial vertex at time $T-1$ when starting with $k_0$ vertices})\cdot \\ & &
\mathbb{P}(\text{being at initial vertex at time $1$ when starting with $k_0+T-1$ }) \\ &+&
\mathbb{P}(\text{not being at initial vertex at time $T-1$ when starting with $k_0$ vertices})\cdot \\ & &
\mathbb{P}(\text{being at  initial vertex at time $1$ when starting from another vertex  with $k_0+T-1$ vertices})    
\end{eqnarray*}
\end{small}
Some elementary manipulations then give
\begin{eqnarray*}
Q(T,k_0) &=& Q(T-1,k_0) Q(1,k_0+T-1) + (1-Q(T-1,k_0))\frac{1}{k_0+T-2}\left(1-Q(1,k_0+T-1)\right) \\
&=& Q(T-1,k_0)\left(\frac{1}{k_0+T-1}+\frac{k_0+T-2}{k_0+T-1}e^{\frac{-\lambda (k_0+T-1)}{k_0+T-2}}\right) \\
&+& (1-Q(T-1,k_0))\frac{1}{k_0+T-2}\left(1-\left(\frac{1}{k_0+T-1}+\frac{k_0+T-2}{k_0+T-1} e^{\frac{-\lambda (k_0+T-1)}{k_0+T-2}}\right)\right)\\
&=& Q(T-1,k_0)e^{\frac{-\lambda (k_0+T-1)}{k_0+T-2}}+ \frac{\left(1- e^{\frac{-\lambda (k_0+T-1)}{k_0+T-2}}\right)}{k_0+T-1}.
\end{eqnarray*}
Division of both sides by $\prod_{i=1}^{T}e^{-\lambda \frac{k_0+i-1}{k_0+i-2}}$ yields
\begin{eqnarray*}
\frac{Q(T,k_0)}{\prod_{i=1}^{T}e^{-\lambda \frac{k_0+i-1}{k_0+i-2}}}=
\frac{Q(T-1,k_0)}{\prod_{i=1}^{T-1}e^{-\lambda \frac{k_0+i-1}{k_0+i-2}}}+
\frac{\left(1-e^{-\lambda\frac{k_0+T-1}{k_0+T-2}}\right)}{(k_0+T-1)\prod_{i=1}^{T}e^{-\lambda \frac{k_0+i-1}{k_0+i-2}}}.
\end{eqnarray*}

Summing up these expressions from $1$ to $T$, and multiplication by $\prod_{i=1}^{T}e^{-\lambda \frac{k_0+i-1}{k_0+i-2}}$ conclude the proof. 
\hfill $\square$
 
Next  we will calculate $E_p^{k_0}$ which is defined to be  the expected number of visits to  vertex $v$  by the end of  the first unit time  when  the probability of starting from vertex $v$ equals to $p$ and the  initial number of vertices equals to $k_0$. 

\begin{proposition}
We have 
$
E_p^{k_0}= \frac{\lambda}{k_0}+\frac{1-pk_0}{k_0^2}\left(1-e^{-\lambda\frac{k_0}{k_0-1}}\right).
$
\end{proposition}

\textbf{Proof.} Recall that $p_r$ is the probability of being at vertex $v$ after $r$ steps. We have
\begin{eqnarray*}
E_p^{k_0}&=& \sum_{r=0}^\infty \mathbb{P}(M_1=r) \mathbb{E}[\text{number of visits to vertex } v \text{ in } r \text{ steps} ]
\\
&=& \sum_{r=0}^\infty \mathbb{P}(M_1=r) \left(\sum_{i=1}^r p_i\right) =\sum_{r=0}^\infty \frac{e^{-\lambda}\lambda^r}{r!}\cdot \left(\sum_{i=1}^r p_i\right)= 
\sum_{r=1}^\infty \frac{e^{-\lambda}\lambda^r}{r!}\cdot \left(\sum_{i=1}^r p_i\right)\\
&=& \sum_{r=1}^\infty \frac{e^{-\lambda}\lambda^r}{r!}\cdot \sum_{i=1}^r \left(\frac{1}{k_0}-\frac{1-pk_0}{k_0}\left(-\frac{1}{k_0-1}\right)^i \right)\\
&=& \sum_{r=1}^\infty \frac{e^{-\lambda}\lambda^r}{r!}\cdot \left( \frac{r}{k_0}+\frac{1-pk_0}{k_0(k_0-1)} \cdot\sum_{i=0}^{r-1} \left(-\frac{1}{k_0-1}\right)^i\right) 
\\&=& \frac{\lambda}{k_0}\sum_{r=1}^\infty \frac{e^{-\lambda}\lambda^{r-1}}{(r-1)!}+ \sum_{r=1}^\infty \frac{e^{-\lambda}\lambda^{r}}{r!}\cdot\frac{1-pk_0}{k_0^2}\cdot\left( 1-\left(-\frac{1}{k_0-1}\right)^r \right)
\\&=&\frac{\lambda}{k_0}+\frac{1-pk_0}{k_0^2}\left(1-e^{-\lambda}+e^{-\lambda}-e^{-\lambda\frac{k_0}{k_0-1}}\right) = \frac{\lambda}{k_0}+\frac{1-pk_0}{k_0^2}\left(1-e^{-\lambda\frac{k_0}{k_0-1}}\right).
\end{eqnarray*}\hfill $\square$

Define $E_p^{k_0}[T]$ to be the expectation of number of visits to vertex $v$ up to time $T$, when the probability of starting from vertex $v$ equals to $p$ and the number of vertices equals to $k_0$ at the beginning.
In particular, $E_p^{k_0}[1] = E_p^{k_0}=  \frac{\lambda}{k_0}+\frac{1-p{k_0}}{{k_0}^2}\left(1-e^{-\lambda\frac{{k_0}}{{k_0}-1}}\right).$
Our main result is on  $E_1^{k_0}[T]$.
\begin{theorem}\label{thm:mainexp}
We have 
$E_1^{{k_0}}[T] = \sum_{i=0}^{T-1} \left(\frac{\lambda}{{k_0}+i}+\frac{1-Q(i,{k_0})({k_0}+i)}{({k_0}+i)^2}\left(1-e^{-\lambda\frac{{k_0}+i}{{k_0}+i-1}}\right)  \right)$,
 where $Q(i,k_0)$'s are as given in Proposition \ref{propn:qtk}.
\end{theorem}

\textbf{Proof.}
 We will sum up the expected values corresponding to  each time interval. The probability of being at vertex $v_1$ at time $i$ is $Q(i,k_0)$. The expected number  of visits to  $v_1$ from time $i$ to $i+1$ is $E_{Q(i,{k_0})}^{{k_0}+i}[1]$. We get
$$
E_1^{{k_0}}[T] = \sum_{i=0}^{T-1}\left(E_{Q(i,{k_0})}^{{k_0}+i}[1] \right) 
= \sum_{i=0}^{T-1} \left(\frac{\lambda}{{k_0}+i}+\frac{1-Q(i,{k_0})({k_0}+i)}{({k_0}+i)^2}\left(1-e^{-\lambda\frac{{k_0}+i}{{k_0}+i-1}}\right)  \right). 
$$
\vspace{-0.1in}
 \hfill $\square$

\vspace{0.1in} 
 \textbf{Acknowledgement.} We would like to thank Nesin Mathematics Village for their kind hospitality where part of this work was done.  The second author is supported partially by BAP grant 20B06P1. We would also like to thank the anonymous referee whose suggestions and
corrections improved the paper significantly.

\end{document}